\documentclass{amsart}[12pt]
\usepackage{latexsym}
\usepackage{amssymb}
\usepackage{amsmath}
\usepackage{amsfonts}
\usepackage{graphicx}
\usepackage{algorithm}
\usepackage{algorithmic}

\global\long\def\Ball{\mathrm{Ball}}
\global\long\def\lspan{\mathrm{span}}
\newcommand{\ignore}[1]{{}}

\def\R{\mathbb{R}}

\newcommand{\bmat}{\left[\begin{matrix}}
\newcommand{\emat}{\end{matrix}\right]}
\usepackage[latin1]{inputenc} 
\usepackage{amsthm}
\usepackage{amsxtra}
\usepackage{amscd}
\usepackage{setspace}   
\usepackage{mathrsfs}   

\newtheorem{theorem}{Theorem}

\numberwithin{theorem}{section}
\numberwithin{equation}{section}

\theoremstyle{remark}

\theoremstyle{definition}

\newcommand{\Z}{\mathbb{Z}}

\renewcommand{\vec}[1]{\mathbf{#1}}

\global\long\def\mQ{\mathbb{Q}}
\global\long\def\mR{\mathbb{R}}

\global\long\def\mZ{\mathbb{Z}}

\global\long\def\vol{\mathrm{vol}}

\global\long\def\round#1{\lfloor#1\rceil}

\global\long\def\norm#1{\left\Vert #1\right\Vert }
\global\long\def\norms#1{\| #1\| }

\newcommand{\vt}{\mathbf{t}}
\newcommand{\vu}{\mathbf{u}}

\usepackage{enumerate}
\usepackage{xcolor}
\usepackage{hyperref}

    \title{On counterexamples to the Mertens conjecture}
    \author[S. Kim]{Seungki Kim}
    \address{Department of Mathematical Sciences, University of Cincinnati, 4199 French Hall West, 2815 Commons Way, Cincinnati, OH 45221-0025, United States}
    \email{seungki.math@gmail.com}
\author[P. Q. Nguyen] {Phong Q. Nguyen}
\address{
Department of Computer Science, Ecole Normale Sup\'erieure, 45 rue d'Ulm, 75005, Paris, France}
\email{phong.nguyen@inria.fr}
 
\begin{document}
\maketitle

\begin{abstract}
We use state-of-art lattice algorithms to improve the upper bound on the lowest counterexample to the Mertens conjecture to $\approx \exp(1.96 \times 10^{19})$, which is significantly below the conjectured value of $\approx \exp(5.15 \times 10^{23})$ by Kotnik and van de Lune \cite{KL04}.
\end{abstract}

\section{Introduction}

The Mertens conjecture \cite{Mertens}, dating back to 1897, is a statement about the growth rate of the Mertens function
\begin{equation*}
M(x) := \sum_{1 \leq n \leq x} \mu(n),
\end{equation*}
where $\mu(n)$ is the M\"obius function
\begin{equation*}
\mu(n) = \begin{cases} (-1)^k & \mbox{if $n$ is squarefree, and has $k$ distinct prime factors,} \\ 0 & \mbox{if $n$ is not squarefree.} \end{cases}
\end{equation*}
The size of $M(x)$ is of interest in number theory, since it is closely related to the size of the real parts of the zeroes of the Riemann zeta function $\zeta(s)$. For example, a short argument (See e.g. \cite[Sec. 2]{OR85}) shows that if $M(x) = O(x^\theta)$, then the Riemann zeta function $\zeta(s)$ has no zeroes on the half-plane $\mathrm{Re}\,s > \theta$. For $\theta = 1/2 + \varepsilon$ for arbitrarily small $\varepsilon>0$, the latter statement is the famous Riemann hypothesis. The Mertens conjecture is a much bolder claim that
\begin{equation}\label{eq:mertens}
|M(x)| < x^{1/2} \mbox{ for all } x>1.
\end{equation}

It took nearly a century for this conjecture to be disproved, by Odlyzko and te Riele \cite{OR85} in 1985. Their argument consisted of certain insights from classical analytic number theory, and, perhaps surprisingly, the use of a lattice reduction algorithm. No alternative (dis)proof that does not rely on lattice reduction is known to this day.

A natural follow-up is to ask about the size of the smallest counterexample $\mathfrak x$ to the Mertens conjecture \eqref{eq:mertens}. This is also related to the estimation of the growth rate of $M(x)$, for which several different conjectures exist. For example, the experimental work by Kotnik and van de Lune \cite{KL04} suggests that
\begin{equation}\label{eq:growth}
|M(x)x^{-1/2}| \approx 1/2 \cdot \sqrt{\log\log\log x}
\end{equation}
along the local extrema of $M(x)/x^{1/2}$, from which they derive the conjecture that 
\begin{equation} \label{eq:conj}
\mathfrak x \approx \exp(5.15 \times 10^{23}),
\end{equation}
since $1/2 \cdot \sqrt{\log\log(5.15 \times 10^{23})} \approx 1$.

It is also possible to give rigorous upper bounds on $\mathfrak x$, thanks to Pintz \cite{Pintz}. The original data of Odlyzko and te Riele \cite{OR85}, under the theorem of Pintz \cite{Pintz}, translates to the statement that $\mathfrak x < \exp(3.21 \times 10^{64})$. The later more extensive set of experiments by Kotnik and te Riele \cite{KR06} led to the improvement $\mathfrak x < \exp(1.59 \times 10^{40})$. Saouter and te Riele \cite{SR14} refined the estimates given in \cite{Pintz}, and also ran more experiments of the same kind, which resulted in $\mathfrak x < \exp(1.004 \times 10^{33})$.

\cite{OR85} states that their achievement was possible thanks to a then-breakthrough in lattice reduction, the LLL algorithm due to Lenstra, Lenstra, and Lov\'asz \cite{LLL}. \cite{KR06} and \cite{SR14} also used LLL, and so did the relatively recent works such as Hurst \cite[Theorem 6.1]{Hurst}. However, over the last decade, there has been a huge amount of progress in lattice reduction techniques --- largely motivated by the emergence of post-quantum cryptography based on computationally hard lattice problems. Noticing this, one would naturally be inclined to apply them to the context of the Mertens conjecture. Kim and Rozmarynowycz \cite{KR23} was the first to point this out and to partially implement this idea, by simply replacing LLL with the BKZ algorithm (\cite{SE94}, \cite{CN11}). As a result, they obtained a further improvement that $\mathfrak x < \exp(1.017 \times 10^{29})$. \cite{KR23} was, however, still far from taking the full advantage of the recent advances in lattice problems.

In this work, we employ state-of-art lattice point enumeration techniques~\cite{LN13,GNR10} in search of a tighter upper bound on the lowest counterexample to the Mertens conjecture \eqref{eq:mertens}. As explained in Section 3 below, this is a much more natural and efficient strategy than running lattice reduction hundreds to thousands of times, which has been the method chosen by all the previous works mentioned above. We also give a careful consideration to the family of the lattices under question, since such an understanding is important for both predicting the outcome and improving the performance of the algorithm. As a result, we improve the choice of the lattice made by \cite{OR85} in a couple of ways, for the first time in the past 38 years. Moreover, we notice that our lattices are extremely orthogonal, with one unusually short vector; we adapt the enumeration algorithm accordingly to exploit this special shape, effectively speeding up our search by a factor of a few millions.

Thereby we were able to find a number of data points that beat the previous record and even the conjecture \eqref{eq:conj} by several orders of magnitude. The best we found, which took us barely half a day on a single core, implies that
\begin{equation} \label{eq:result}
\mathfrak x < \exp(1.96 \times 10^{19}),
\end{equation}
which is below \eqref{eq:conj} by a factor of $\approx 26276$ in the exponent. This would substantially impact the credibility of \eqref{eq:growth}. A more extensive application of the methods of the present paper may help search for an alternative estimate on the growth rate of $M(x)$ --- see Section 4.2 below.

There exist also numerous other applications of the computational aspect of lattices to number theory than the Mertens conjecture, see e.g. Simon \cite{Simon}. We hope that our work helps inform the community of the recent advances in lattice computations, and encourages revisiting some of the old problems with the new arsenal.

\subsection*{Organization}

In Section 2, we briefly review the method of Odlyzko-te Riele \cite{OR85} and other works in the literature, in order to help the reader understand and put our work into perspective. In Section 3, we describe our experiment in detail. We conclude the paper in Section 4 with a discussion on further research directions to which the methods introduced herein may be applied.

\subsection*{Reproducibility}
The present work involves significant data.
Data files and/or source codes allowing to reproduce the data are available on \url{https://zenodo.org/records/10775723}.

\subsection*{Acknowledgments}

This project has received funding from the European Research Council (ERC) under the European Union's Horizon 2020 research and innovation programme (grant agreement No 885394).
In addition, S.K. was supported by NSF grant CNS-2034176.
We thank the anonymous referees for helpful comments and suggestions.

\section{A review of the previous works}

\subsection{The original argument of \cite{OR85}}

Resorting to proof by contradiction, let us assume the truth of the Mertens conjecture. As explained in \cite[Sec. 2]{OR85}, this implies in particular that all nontrivial zeroes $\rho$ of $\zeta(s)$ satisfy $\mathrm{Re}\, \rho = 1/2$ (the Riemann hypothesis), and that they are all simple zeroes. For each such $\rho$ with $\mathrm{Im}\,\rho > 0$, let us write $\gamma := \mathrm{Im}\,\rho$, $\alpha := |\rho\zeta'(\rho)|^{-1}$, $\psi := \mathrm{arg}(\rho\zeta'(\rho))$; conversely, whenever we write $\gamma$, $\alpha$, or $\psi$, we are referring to the corresponding zero $\rho$ of $\zeta(s)$. Then, for a certain increasing sequence $\{T_n\}$ with $n \leq T_n \leq n+1$, it holds that (see e.g. \cite[(3) and (4)]{KR06} or \cite[Sec. 2 and 3]{KL04})
\begin{equation}\label{eq:q}
q(x) := \frac{M(x)}{x^{1/2}} = 2\lim_{n \rightarrow \infty} \sum_{0 < \gamma < T_n} \alpha\cos(\gamma y - \psi) + O(x^{-1/2}),
\end{equation}
where we write $y := \log x$ for short. Recall that the Mertens conjecture states that $|q(x)| < 1$ for $x>1$.

\eqref{eq:q} suggests one possible strategy for disproving \eqref{eq:mertens}: for a large $N>0$, find $y$ such that the sum
\begin{equation}\label{eq:qN}
q_N(x) := 2\sum_{0 < \gamma < N} \alpha \cos(\gamma y - \psi)
\end{equation}
is large. This can be interpreted as a problem in simultaneous Diophantine approximation, as follows. Let us denote by $|a|_{2\pi}$ the representative of $a \mbox{ (mod  $2\pi$)}$ in $(-\pi,\pi]$. If we can find $y$ such that all $|\gamma y - \psi|_{2\pi}$ are small, then we can expect that
\begin{align*}
2\sum_{0 < \gamma < N} \alpha \cos(\gamma y - \psi) &\approx \sum_{0 < \gamma < N} \alpha(2-|\gamma y - \psi|^2_{2\pi}) \\
&\approx \sum_{0 < \gamma < N} 2\alpha.
\end{align*}
Now it is known that the last sum diverges as $N \rightarrow \infty$. Hence, for a sufficiently large $N$, if we can indeed find such $y$ and thereby not lose too much in the estimates above, we can hope to be able to demonstrate that \eqref{eq:q} is greater than $1$, and --- at least in principle --- even arbitrarily large. This is precisely the approach taken by \cite{OR85}.

(Or alternatively, by finding $y$ such that all $|\gamma y - \psi - \pi|_{2\pi}$ are small, we can try to show $q(x)$ can be large in the negative direction.)

\cite{OR85} converts the problem of the (weighted inhomogeneous simultaneous) Diophantine approximation to the problem of reducing a certain lattice. They consider the lattice in $\R^{N+2}$, say $L_0$, consisting of the integer linear combinations of the rows of
\begin{equation} \label{eq:otrmat}
\begin{pmatrix}
\lfloor 2\pi\sqrt{\alpha_1}2^\nu \rfloor & 0 & \cdots &0 & 0 & 0 \\

0 &  \lfloor 2\pi\sqrt{\alpha_2}2^\nu \rfloor & \cdots & 0 & 0 & 0 \\

\vdots & \vdots & \ddots & \vdots & \vdots & \vdots \\

0 & 0 & \cdots & \lfloor 2\pi\sqrt{\alpha_N}2^\nu \rfloor & 0 & 0 \\

-\lfloor \sqrt{\alpha_1}\psi_12^\nu \rfloor & -\lfloor \sqrt{\alpha_2}\psi_22^\nu \rfloor &  \cdots & -\lfloor \sqrt{\alpha_N}\psi_N2^\nu \rfloor & 2^\nu N^4 & 0 \\

\lfloor \sqrt{\alpha_1}\gamma_12^{\nu-10} \rfloor & \lfloor \sqrt{\alpha_2}\gamma_22^{\nu-10} \rfloor & \cdots & \lfloor \sqrt{\alpha_N}\gamma_N2^{\nu-10} \rfloor & 0 & 1 
\end{pmatrix},
\end{equation}
where $\alpha_i$'s are the $\alpha$'s ordered in descending order so that $\alpha_1 > \alpha_2 > \ldots$, and the $\psi_i$, etc., are those corresponding to the zero $\rho_i$ of $\zeta(s)$ associated to $\alpha_i$; and the role of $2^\nu$ is to approximate the entries of \eqref{eq:otrmat} to $\nu$ most significant base $2$ digits.
Given a basis of a lattice such as this one, lattice reduction algorithms such as LLL compute another basis of the same lattice consisting of vectors of reasonably short length (depending on the strength of the algorithm), called a \emph{reduced basis} of that lattice.

Let us consider a reduced basis of $L_0$. We claim that it must contain a vector of the form
\begin{align}
&(p_1\lfloor 2\pi\sqrt{\alpha_1}2^\nu \rfloor + z \lfloor \sqrt{\alpha_1}\gamma_12^{\nu-10} \rfloor -\lfloor \sqrt{\alpha_1}\psi_12^\nu \rfloor, \ldots \label{eq:otrvr} \\ 
&\ldots, p_N\lfloor 2\pi\sqrt{\alpha_N}2^\nu \rfloor + z \lfloor \sqrt{\alpha_N}\gamma_N2^{\nu-10} \rfloor -\lfloor \sqrt{\alpha_N}\psi_N2^\nu \rfloor, \pm2^\nu N^4, z) \notag
\end{align}
for some integers $p_1, \ldots, p_N$ and $z$. There certainly must be a vector whose penultimate entry is nonempty, since a reduced basis is, in particular, a basis. But then, since $2^{\nu}N^4$ is huge compared to the rest of the entries of \eqref{eq:otrmat}, one can argue from the performance guarantee of LLL \cite[Prop. 1.12]{LLL} that the penultimate entry must be as small as possible.

Now set $y = \pm2^{10}z$, the sign being that of $2^{\nu}N^4$ in \eqref{eq:otrvr}. Then it can be seen, from the ``size-reducedness'' property of a reduced basis \cite[(1.4)]{LLL}, that each of the first $N$ entries of \eqref{eq:otrvr} are approximately equal to $2^\nu\sqrt{\alpha} \cdot |\gamma_iy - \psi|_{2\pi}$.

One may still be (rightfully) curious as to how this would result in a good Diophantine approximation $y$. As \cite[Sec. 3]{OR85} writes, some part of it is a miracle: LLL is famously known to perform much better than the theoretical guarantee given by \cite[Prop. 1.12]{LLL}. However, as pointed out in \cite{KR23}, much of it can be adequately explained in the language of lattice problems. What the reduction of \eqref{eq:otrmat} really achieves is the resolution of the \emph{approximate closest vector problem} (aCVP) of finding a point on the lattice in $\R^{N+1}$ generated by the rows of
\begin{equation} \label{eq:otrmat'}
\begin{pmatrix}
\lfloor 2\pi\sqrt{\alpha_1}2^\nu \rfloor & 0 & \cdots &0 & 0 \\

0 &  \lfloor 2\pi\sqrt{\alpha_2}2^\nu \rfloor & \cdots & 0 & 0 \\

\vdots & \vdots & \ddots & \vdots & \vdots \\

0 & 0 & \cdots & \lfloor 2\pi\sqrt{\alpha_N}2^\nu \rfloor & 0 \\

\lfloor \sqrt{\alpha_1}\gamma_12^{\nu-10} \rfloor & \lfloor \sqrt{\alpha_2}\gamma_22^{\nu-10} \rfloor & \cdots & \lfloor \sqrt{\alpha_N}\gamma_N2^{\nu-10} \rfloor & 1 
\end{pmatrix}
\end{equation}
that is reasonably close to the ``target'' vector
\begin{equation*}
(-\lfloor \sqrt{\alpha_1}\psi_12^\nu \rfloor, -\lfloor \sqrt{\alpha_2}\psi_22^\nu \rfloor,  \cdots, -\lfloor \sqrt{\alpha_N}\psi_N2^\nu \rfloor, 0) \in \R^{N+1}
\end{equation*}
via \emph{Babai's nearest plane algorithm} \cite{Babai}, one of the standard approaches to aCVP to this day, and implicitly used as a subroutine inside LLL itself. This interpretation allows one to predict the outcome heuristically, that matches the actual output rather well --- see \cite[Sec. 2]{KR23} for details.

\subsection{Results on the smallest counterexample}

As in the introduction, let us continue to denote by $\mathfrak x$ the smallest real number for which \eqref{eq:mertens} does not hold.
Giving a rigorous upper bound on $\mathfrak x$ became possible thanks to the following result of Pintz \cite{Pintz}.

\begin{theorem}[Pintz \cite{Pintz}] \label{thm:pintz}
Let
\begin{equation} \label{eq:hp}
h_P(y) := 2\sum_{\gamma < 14000} \alpha\exp(-1.5 \cdot 10^{-6}\gamma^2)\cos(\gamma y - \psi).
\end{equation}
If there exists $y \in [e^7, e^{50000}]$ with $|h_P(y)| > 1+e^{-40}$, then $\mathfrak x < \exp(y + \sqrt{y})$.
\end{theorem}

For the value of $y$ found in \cite{OR85}, Theorem \ref{thm:pintz} implies that $\mathfrak x < \exp(3.21 \times 10^{64})$. Later, \cite{KR06} repeated the experiment of \cite{OR85}, that we described in the previous section, over a broader range of parameters $N$ and $\nu$, and found a lower working value of $y$, which corresponds to the improved bound $\mathfrak x < \exp(1.59 \times 10^{40})$.

\cite{SR14} improved Theorem \ref{thm:pintz} by refining Pintz's estimates on certain contour integrals involving the zeta function. Consequently, they were essentially able to replace $h_P$ in Theorem \ref{thm:pintz} by
\begin{equation} \label{eq:hstr}
h_{StR}(y) := 2\sum_{\gamma < 74000} \alpha\exp(-3 \cdot 10^{-9}\gamma^2)\cos(\gamma y - \psi).
\end{equation}
As can be seen by comparing \eqref{eq:hp} and \eqref{eq:hstr}, $|h_{StR}|$ tends to be somewhat larger than $|h_P|$, so that some of the ``near-misses,'' i.e. those $y$ for which $|h_P(y)| < 1$ but very close, may satisfy $|h_{StR}(y)| > 1$ and become valid ``hits.'' For practical values of $y$, $|h_{StR}(y)| > 1 + 6 \cdot 10^{-8}$ implies that $\mathfrak x < \exp(y + \sqrt{y})$.
With this, and some extra search for the candidate values, \cite{SR14} attained $\mathfrak x < \exp(1.004 \times 10^{33})$.

Recently, \cite{KR23} essentially repeated the experiment of \cite{KR06}, except for the following few tweaks:
\begin{enumerate}[(i)]
\item They replaced LLL with the more powerful BKZ, which leads to a much better solution to the aCVP problem.
\item They ran tens of thousands of trials, and for each trial, they perturbed the basis \eqref{eq:otrmat} hoping that the randomization effect would help find a lower value of $y$ for which Theorem \ref{thm:pintz} is applicable.
\end{enumerate}
These were possible --- with only the computational power of a personal laptop --- thanks to the substantial advances in lattice reduction over the last decade.
They succeeded in finding the value
\begin{equation*}
y = 1017256208\ 7569945816\ 8018857216.806640625 \mbox{ with } h_P(y) = 1.0034372\ldots,
\end{equation*}
which implies $\mathfrak x < \exp(1.017 \times 10^{29})$, the best record to this date.


\section{Our experiment}

\subsection{Lattice point enumeration}

The method of lattice reduction for finding $y$ satisfying the conditions of Theorem \ref{thm:pintz} has a few limitations. To maximize \eqref{eq:hp} or \eqref{eq:hstr}, it makes sense to account for as many summands as possible, that is, to take $N$ as large as possible. However, high-quality lattice reduction becomes extremely time-consuming for $N \geq 100$, so the trial-and-error strategy of \cite{KR23} would take too much computational resource: non-trivial lattice tasks such as finding a closest lattice point typically run in time exponential in $N$.
Furthermore, maximizing \eqref{eq:hp} or \eqref{eq:hstr} is not exactly a lattice problem: 
it is possible that the optimal solution does not arise from a particularly close lattice point.
The terms on the right-hand side of \eqref{eq:hp} or \eqref{eq:hstr} that are left out in the construction \eqref{eq:otrmat} turn out to be large enough to affect the outcome either in or against our favor, and they are essentially the matter of a coin toss.

The technique of lattice point enumeration provides a natural solution to this dilemma.
When the lattice dimension is moderate, it is possible to enumerate all the lattice points within any small ball,
or if the ball is larger, to enumerate many lattice points inside the ball.
To do this, a standard algorithm is enumeration with pruning~\cite{LN13,GNR10},
building upon work by Pohst~\cite{P81}, Kannan~\cite{K83}, and Schnorr-Euchner~\cite{SE94}.
Gama, Nguyen and Regev~\cite{GNR10} showed how to speed up rigorously the Schnorr-Euchner~\cite{SE94}
enumeration of lattice points within a zero-centered ball, which consists of a depth-first search (DFS) of a carefully constructed tree. This technique was later adapted to any ball  by Liu and Nguyen~\cite{LN13}.

Let $L$ be a full-rank lattice in $\mR^n$.
Given a target $\vec{t} \in \mQ^n$, a basis $B=(\vec{b}_1,\dots,\vec{b}_n)$ of $L$ and a radius $R>0$, 
enumeration~\cite{P81,K83} outputs $L \cap S$ where $S=\Ball_{n}(\vec{t},R)$:.
It performs a recursive search using projections, to reduce the dimension of the lattice:
if $\|\vec{v}\| \le R$, then $\|\pi_k(\vec{v})\| \le R$ for all $1 \le k \le n$,
where $\pi_{k}$ denotes the orthogonal projection on $\lspan(\vec{b}_{1},\dots,\vec{b}_{k-1})^{\perp}$.
One can easily enumerate  $\pi_n(L) \cap S$.
And if one enumerates  $\pi_{k+1}(L) \cap S$ for some $k \ge 1$,
one derives  $\pi_{k}(L) \cap S$ by enumerating
the intersection of a one-dimensional lattice with a suitable ball,
for each point in $\pi_{k+1}(L) \cap S$.
Concretely, it can be viewed as a depth-first search of the following enumeration tree $\mathcal{T}$:
the nodes at depth $n+1-k$ are the points of $\pi_{k}(L) \cap S$.
The  running-time of enumeration depends on $R$ and $B$, but is typically super-exponential in $n$, even if $L\cap S$ is small.

Pruned enumeration~\cite{GNR10,SE94} uses a pruning set $P \subseteq \mR^n$,
and outputs  $L \cap (\vec{t}+P)$.
The advantage is that for suitable choices of $P$,
enumerating  $L \cap (\vec{t}+P)$ is much cheaper than $L \cap S$, yet under mild heuristics, $L \cap (\vec{t}+P)$ is expected to cover most of  $L \cap S$.
The pruning set $P$ should be viewed as a random variable:  it depends on the choice of basis $B$.
In~\cite{GNR10}, $P$ is defined by a function $f: \{1,\dots,n\} \rightarrow [0,1]$, a radius $R>0$ and a lattice basis $B=(\vec{b}_1,\dots,\vec{b}_n)$
as follows:
\begin{align} \label{eqn:defPf}
 P_f(B,R) = \{ \vec{x} \in \mR^n \,\,\text{s.t.}\,\, \|\pi_{n+1-i}(\vec{x})\| \le f(i) R \,\,\text{for all}\,\, 1 \le i \le n \},
 \end{align}
 This form of pruning is known as cylinder pruning,
because $P_f(B,R)$ is an intersection of cylinders: each inequality $\|\pi_{n+1-i}(\vec{x})\| \le f(i) R$ defines a cylinder.
And \cite{GNR10} provides an algorithm which, given as input $(B,R,f)$, decides if $L \cap P_f(B,R)$ is non-empty, and if so, outputs one element of $L \cap P_f(B,R)$.
It is easy to modify this algorithm to tackle more needs: for instance, \cite{LN13} extended it to  $L \cap (\vec{t}+P_f(B,R))$where $\vec{t}$ is an additional input.
In our case, we are interested in enumerating the whole $L \cap (\vec{t}+P_f(B,R))$:
specifically, we implemented Alg.~\ref{alg:prunedEnum}, which is a slight variant of \cite{LN13}.
\begin{algorithm}
\caption{Pruned Enumeration for BDD of unbalanced lattices \label{alg:prunedEnum} (slight variant version of~\cite{LN13,GNR10})}

\begin{algorithmic}[1]
\REQUIRE A reduced basis $B=(\vec{b}_{1},\dots,\vec{b}_{m})$ of a lattice $L$, a target vector $\vec{t}=\sum_{i=1}^m t_i \vec{b}_i$,
a bounding function $R_{1}^{2}\leq\dots\leq R_{m}^{2}$,
the Gram-Schmidt matrix $\mu = (\langle \vec{b}_i,\vec{b}_{j}^{*} \rangle / \norms{\vec{b}_{j}^{*}}^2)$
and the (squared) norms $\norms{\vec{b}_{1}^{*}}^2, \ldots, \norms{\vec{b}_{m}^{*}}^2$,
where the $\vec{b}_{j}^*$'s are the Gram-Schmidt orthogonalization of the $\vec{b}_{j}$'s.
\ENSURE The basis coefficients of all $\vec{v} \in L$ such that
the projections of $\vec{v}-\vec{t}$ have norms less than the $R_i$'s,
{\em i.e.} $\norm{\pi_{m+1-k}(\vec{v}-\vec{t})}\leq R_k$ for all $1 \le k \le m$,
and there is no $x_1 \in \mZ$ such that $\|\vec{v}+ x\vec{b}_1-\vec{t} \| <  \|\vec{v}-\vec{t} \|$.
The latter constraint avoids returning too many vectors, when $\vec{b}_1$ is much shorter than $\det(L)^{1/m}$.

\STATE \label{line:optimization1}$\sigma\leftarrow(0)_{(m+1)\times m}$; $r_0 = 0; r_1 = 1; \cdots ; r_m = m$;$\rho_{m+1}=0$

\STATE \textbf{for} $k$ = $m$ \textbf{downto}
$1$
\STATE
\ \ \ \textbf{for} $i$ = $m$ \textbf{downto}
$k+1$ \textbf{do} $\sigma_{i,k}\leftarrow\sigma_{i+1,k}+(t_{i}-v_{i})\mu_{i,k}$
\textbf{endfor}
\STATE \ \ \ $c_{k} \leftarrow t_k+\sigma_{k+1,k}$ \emph{ // $c_k \leftarrow t_{k}+ \sum_{i=k+1}^m (t_i-v_i) \mu_{i,k}$,  \emph{ centers}}
\STATE  \ \ \ $v_{k}  \leftarrow \round{c_{k}}$\emph{ // current combination};
\STATE \ \ \ $w_{k}=1$\emph{ // jumps};
\STATE \ \ \ $\rho_{k}=\rho_{k+1}+(c_{k}-v_{k})^{2}\cdot\norms{\vec{b}_{k}^{*}}^2$
\STATE \textbf{endfor}
\STATE $k=1$;
\WHILE{true}
 \STATE $\rho_{k}=\rho_{k+1}+(c_{k}-v_{k})^{2}\cdot\norms{\vec{b}_{k}^{*}}^2$
   \emph{ // compute norm squared of current node}
 \IF{\label{line:boundingfunc}$\rho_{k}\leq R_{m+1-k}^{2}$\emph{ }(\emph{we are
below the bound})}
   \IF{$k=1$}
     \STATE \textbf{return} $(v_1,\ldots,v_m)$; (\emph{solution found; so we're going up the tree})
      \STATE $k\leftarrow k+1$ \emph{// going up the tree}
   \IF{{$k = m+1$}}
     \STATE \textbf{return} $\emptyset$ (\em{there is no solution})
   \ENDIF
   \STATE \label{line:optimization5} $r_{k-1} \leftarrow k$ \emph{// since $v_k$ is about to change,
   indicate that $(i,j)$ for $j < k$ and $i \le k$ are not synchronized}
   \STATE \emph{// update $v_k$}
     \STATE \textbf{if} $v_{k}>c_{k}$
     \textbf{then} $v_{k}\leftarrow v_{k}-w_{k}$
     \textbf{else} $v_{k}\leftarrow v_{k}+w_{k}$
     \STATE $w_{k}\leftarrow w_{k}+1$
   \ELSE
     \STATE $k\leftarrow k-1$ \emph{// going down the tree}
     \STATE \label{line:optimization2} $r_{k-1} \leftarrow \max(r_{k-1},r_k)$ \emph{// to maintain the invariant for $j<k$}
   \STATE \label{line:optimization3}\textbf{for} $i$ = $r_{k}$ \textbf{downto}
$k+1$ \textbf{do} $\sigma_{i,k}\leftarrow\sigma_{i+1,k}+(t_{i}-v_{i})\mu_{i,k}$
\textbf{endfor}
   \STATE \label{line:optimization4}$c_{k} \leftarrow t_k+\sigma_{k+1,k}$ \emph{ // $c_k \leftarrow t_{k}+ \sum_{i=k+1}^m (t_i-v_i) \mu_{i,k}$}
   \STATE  $v_{k}  \leftarrow \round{c_{k}}$;
$w_{k}=1$

   \ENDIF
 \ELSE
   \STATE $k\leftarrow k+1$ \emph{// going up the tree}
   \IF{{$k = m+1$}}
     \STATE \textbf{return} $\emptyset$ (\em{there is no solution})
   \ENDIF
   \STATE \label{line:optimization5} $r_{k-1} \leftarrow k$ \emph{// since $v_k$ is about to change,
   indicate that $(i,j)$ for $j < k$ and $i \le k$ are not synchronized}
   \STATE \emph{// update $v_k$}
     \STATE \textbf{if} $v_{k}>c_{k}$
     \textbf{then} $v_{k}\leftarrow v_{k}-w_{k}$
     \textbf{else} $v_{k}\leftarrow v_{k}+w_{k}$
     \STATE $w_{k}\leftarrow w_{k}+1$
 \ENDIF
\ENDWHILE
\end{algorithmic}
\end{algorithm}

Gama {\em et al.}~\cite{GNR10} showed how to efficiently compute tight lower and upper bounds
for $\vol(P_f(B,R))$  and also estimated of the cost of enumerating $L \cap S \cap P_f(B,R)$,
 using the Gaussian heuristic (that for a lattice $L$ and a set $S$, $|L \cap S| \approx \vol\, S/\det L$) on projected lattices $\pi_i(L)$: these estimates are usually accurate in practice,
 and they can also be used in the CVP case~\cite{LN13}.

In our experiments, we used a function $f$ close to linear.
It is well-known that if $\vec{x}$ denotes the projection of a random unit vector of $\mR^n$ onto an $i$-dimensional subspace,
then $\|\vec{x}\|^2$ has distribution Beta$(i/2,(n-i)/2)$.
Accordingly, we took $f^2(i) = \mu_i + 2 \sigma_i$ where $\mu_i$ and $\sigma_i$ are respectively the expectation and the standard deviation
of the Beta$(i/2,(n-i)/2)$ distribution. From the analysis of~\cite{GNR10}, 
for this choice of $f$,  $L \cap (\vec{t}+P_f(B,R))$ should cover most of $L \cap \Ball_{n}(\vec{t},R)$.

\subsection{Choice of lattice}

With optimizing \eqref{eq:hp} or \eqref{eq:hstr} in mind, for each nontrivial zero $\rho$ of $\zeta(s)$ with $\mathrm{Im}\,\rho > 0$,
we define $\alpha^* = \alpha\exp(-1.5 \cdot 10^{-6}\gamma^2)$ if we want to find large values of $|h_P|$,
or $\alpha^* = \alpha\exp(-3 \cdot 10^{-9}\gamma^2)$ if we want to find large values of of $|h_{StR}|$. We define the corresponding $\gamma$ and $\psi$ as earlier. We index the $\rho$'s and other variables accordingly so that $\alpha_1^* > \alpha_2^* > \ldots$.

For parameters $\nu, \nu_y, \nu_t$, we consider the integer lattice $L$ spanned by the rows of the $(N+1) \times (N+1)$ matrix
\begin{equation} \label{eq:ourmat}
\begin{pmatrix}
\lfloor  \sqrt{\alpha^*_1} \gamma_1 2^{\nu_y} \rfloor &  \lfloor \sqrt{\alpha^*_2} \gamma_2 2^{\nu_y} \rfloor   & \dots &   \lfloor  \sqrt{\alpha^*_N} \gamma_N 2^{\nu_y} \rfloor  &  2^{\nu_t} \\
 \lfloor    \sqrt{\alpha^*_1} 2\pi   2^{\nu}   \rfloor & 0           & \dots    & 0  & 0\\
0           & \ddots   & \ddots & \vdots     & \vdots \\
\vdots   &        \ddots                               &                     & 0    & 0  \\
0         & \dots.                        & 0  &\lfloor   \sqrt{\alpha^*_N} 2\pi  2^{\nu}  \rfloor & 0
\end{pmatrix}.
\end{equation}
We used $\alpha^*$ instead of $\alpha$ as in \eqref{eq:otrmat'}, since we are trying to optimize $h_P$ \eqref{eq:hp} or $h_{StR}$ \eqref{eq:hstr}, not $q_N$ \eqref{eq:qN} as in \cite{OR85}. The new parameters $\nu_y$ and $\nu_t$ will be helpful in further analyzing and improving the outcome of the experiment, as explained in the next section. Our ``target vector'' is
\begin{equation*}
\vt = \left(
\lfloor \sqrt{\alpha^*_1} \psi_1 2^{\nu} \rfloor, \lfloor  \sqrt{\alpha^*_2} \psi_2  2^{\nu}  \rfloor, \cdots,  \lfloor \sqrt{\alpha^*_N} \psi_N 2^{\nu} \rfloor,  0  \\
\right)
\end{equation*}
for positive values of $h_P$ (and
\begin{equation*}
\vt' =
\left(
\lfloor \sqrt{\alpha^*_1} (\psi_1+\pi) 2^{\nu} \rfloor, \lfloor  \sqrt{\alpha^*_2} (\psi_2+\pi))  2^{\nu}  \rfloor, \cdots,  \lfloor \sqrt{\alpha^*_N} (\psi_N+\pi) 2^{\nu} \rfloor, 0  \\
\right)
\end{equation*}
for negative values of $h_P$). For a parameter $\gamma > 0$, we apply lattice enumeration to find all points $\vu \in L$, or equivalently all $x \in \Z$, such that
\begin{align}
\vu-\vt = \Big( 
& \left|x \lfloor \sqrt{\alpha_1^*} \gamma_1 2^{\nu_y} \rfloor - \lfloor 2^\nu  \sqrt{\alpha_1^*} \psi_1 \rfloor\right|_{\lfloor 2^\nu 2\pi \sqrt{\alpha_1^*} \rfloor },
\cdots, \label{eq:u-t} \\
& \cdots, \left|x \lfloor \sqrt{\alpha_N^*} \gamma_N 2^{\nu_y} \rfloor - \lfloor 2^\nu  \sqrt{\alpha_N^*} \psi_N \rfloor\right|_{\lfloor 2^\nu 2\pi \sqrt{\alpha_N^*} \rfloor },
x2^{\nu_t} \Big) \notag
\end{align}
(note the similarity with \eqref{eq:otrvr}) is shorter than
\begin{equation*}
K:= \gamma\sqrt{\frac{N+1}{2\pi e}}\det L^\frac{1}{N+1},
\end{equation*}
where $\gamma \geq 1$ is a parameter to be chosen later, and $\sqrt{\frac{N+1}{2\pi e}}$ is the approximate radius of the ball in $\R^{N+1}$ of unit volume. Hence
\begin{equation} \label{eq:norm1}
 \sum_{i=1}^N \left|x \lfloor \sqrt{\alpha_i^*} \gamma_i 2^{\nu_y} \rfloor - \lfloor 2^\nu  \sqrt{\alpha_i^*} \psi^*_i \rfloor\right|_{\lfloor 2^\nu 2\pi \sqrt{\alpha_i^*} \rfloor }^2
 < K^2,
\end{equation}
and thus we would expect
\begin{equation} \label{eq:norm2}
\sum_{i=1}^N \alpha^*_i |\gamma_i^* y - \psi^*_i |_{2\pi}^2 < \frac{K^2}{2^{2\nu}},
\end{equation}
where $y = x2^{\nu_y - \nu}$.

\begin{figure}[h]
	\includegraphics[height=8cm]{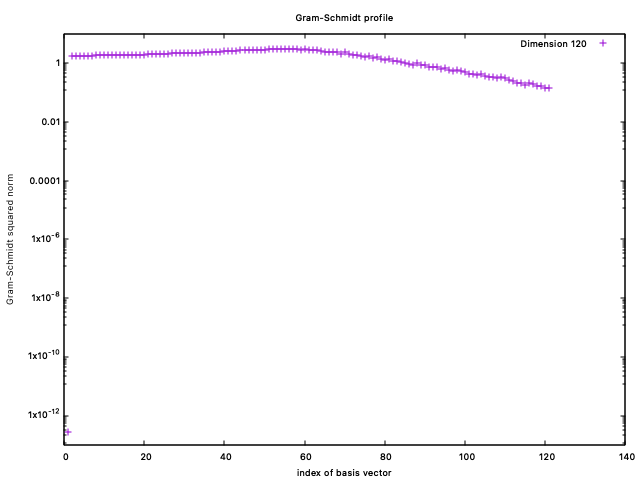}
	 \caption{Profile of a 1-round BKZ-84 reduced basis of the Mertens lattice for $N=120$, $\nu = 130$, $\nu_y = 100$, $\nu_t = 15$, and $\alpha^* = \alpha\exp(-1.5 \cdot 10^{-6}\gamma^2)$. \label{fig:profile}}
\end{figure}

We note that the lattice defined by \eqref{eq:ourmat} is unbalanced when $\nu \gg \nu_y$.
Indeed, in such a case, the $i$-th coefficient of the first row is much smaller than the $i$-th diagonal coefficient and therefore,
the first row has norm much smaller than $\det(L)^{1/(N+1)}$ and is likely to be a shortest vector of $L$:
Fig.~\ref{fig:profile} shows the typical profile of a reduced basis of $L$, which differs from a reduced basis of a random lattice (in the sense of the Haar measure on $\mathrm{PGL}(N+1,\Z)\backslash\mathrm{PGL}(N+1,\R)$).
Here, the first vector is much smaller than $\det(L)^{1/(N+1)}$, and the Gram-Schmidt norms of the reduced basis first increase,
before eventually decreasing geometrically as in a typical reduced basis~\cite{CN11,GNR10}:
this means reduced bases of $L$ are significantly more reduced than a reduced basis of a random lattice,
which makes enumeration faster for the same dimension.
This property must also be taken into account when enumerating lattice points inside a ball.
Indeed, whenever we have found $\vec{u} \in L$ such that $\|\vt-\vec{u}\| \le R$
then there are likely many integers $m \in \mZ$ such that $\|\vt-\vec{u}-m\vec{b}_1\| \le R$, where $\vec{b}_1$ is the top row of \eqref{eq:ourmat},
because $\|\vec{b}_1\|$ is much smaller than $R$ and $\|\vt-\vec{u}\|$.
So if we want to enumerate $L \cap \Ball_{n}(\vec{t},R)$, it is better to only enumerate the points $\vec{u} \in  L \cap \Ball_{n}(\vec{t},R)$
such that there is no nonzero  $m \in \mZ$ such that $\|\vt-\vec{u}-m\vec{b}_1\| < \|\vt-\vec{u}\|$:
this is done by Alg.~\ref{alg:prunedEnum}, which is a slight variant of \cite{LN13}.
In other words, we are actually enumerating over the projection of $L$ onto the hyperplane orthogonal to $\vec{b}_1$.
This phenomenon is inherent to our construction: if $y \approx z$, then $h_P(y) \approx h_P(z)$,
so there is no need to return all the lattice points corresponding to $z \approx y$, whenever we have found a good $y$.

\subsection{Constraints on the parameters}

The above discussion leaves the five parameters $N, \nu, \nu_y, \nu_t, \gamma$ to be determined. In principle, $N$ is the bigger the better, and the only constraint is the amount of computational power available. $\gamma$ controls the number of the candidate points that are to be enumerated, since $|L \cap \Ball_{n}(\vec{t},K)| \approx \gamma^{N+1}$. In our experiments, we chose $N \in \{120, 130, 140 \}$, and $\gamma \in [1,1.28]$.

$\nu, \nu_y, \nu_t$ have a strong influence on the entry sizes of \eqref{eq:u-t}, which helps us choose their values to some extent. We expect $\|\vu-\vt\| \approx K$, since most of the mass of a high-dimensional ball lies away from its center. Therefore, each entry would be of size around
\begin{equation*}
\frac{K}{\sqrt{N+1}}
= \frac{\gamma}{\sqrt{2\pi e}}\det L^{\frac{1}{N+1}}
= \frac{\gamma}{\sqrt{2\pi e}} \cdot 2^\frac{\nu_t}{N+1} 2^\frac{N\nu}{N+1} \prod_{i=1}^N (2\pi\sqrt{\alpha^*_i})^\frac{1}{N+1}.
\end{equation*}
For $N=120$ for instance, we have
\begin{equation*}
\frac{1}{\sqrt{2\pi e}}\prod_{i=1}^N (2\pi\sqrt{\alpha^*_i})^\frac{1}{N+1} \approx 2^{-3.6},
\end{equation*}
from which we can predict
\begin{equation} \label{eq:pred_a}
\sqrt{\alpha^*_i}|\gamma_iy - \psi_i|_{2\pi} \in 2^{-3.6 - \frac{\nu-\nu_t}{N+1}} \cdot [1,\gamma]
\end{equation}
for each $i$. Since we wish this to be small, $\nu - \nu_t$ needs to be at least of comparable size to $N$. By a similar computation, we also obtain the heuristic
\begin{equation} \label{eq:pred_y}
y \in 2^{-3.6 + \nu_y - \nu_t - \frac{\nu-\nu_t}{N+1}} \cdot [1,\gamma].
\end{equation}
Of course, these statements must be taken with a grain of salt, since we are looking at a ball of a relatively small radius, and our lattice has a rather unusual shape, as remarked earlier. Still, they can and do serve as useful guides in practice: in our experiments, the values of $y$ found hardly differed from \eqref{eq:pred_y} by more than a factor of $2^4$, as can be checked from Tables \ref{fig:besthP} and \ref{fig:besthStR} below.

In addition, for the expectation \eqref{eq:norm2} made from \eqref{eq:norm1} to be reasonable, it is necessary that $\nu_t$ be not too small. The reason is that the difference between $x \lfloor \sqrt{\alpha_i^*} \gamma_i 2^{\nu_y} \rfloor$ and $\lfloor x\sqrt{\alpha_i^*} \gamma_i 2^{\nu_y} \rfloor$ can scale with $x$: for instance, $100 \cdot \lfloor 1.99 \rfloor = 100$, whereas $\lfloor 100 \cdot 1.99 \rfloor = 199$. From \eqref{eq:pred_y}, we find $x \approx 2^{-3.6 + \nu - \nu_t - \frac{\nu-\nu_t}{N+1}}$, dividing which by $2^\nu$ gives $2^{-3.6 - \nu_t - \frac{\nu-\nu_t}{N+1}}$. Hence \eqref{eq:norm2} may diverge from our expectation by as much as $N2^{-3.6 - \nu_t - \frac{\nu-\nu_t}{N+1}}$, which may be nontrivial for small values of $\nu_t$. The previously used lattice construction \eqref{eq:otrmat'} was problematic in this respect, since it fixes $\nu_t=0$.

\subsection{Implementation details}

We used three software libraries: Arb~\cite{ARB} for guaranteed interval arithmetic to compute $h_P$ and $h_{StR}$, 
fplll~\cite{FPLLL} for lattice basis reduction (implementations of LLL and BKZ), and NTL~\cite{NTL} with which we implemented
our variant of the pruned enumeration algorithm~\cite{LN13}.

For the values $\rho, \alpha, \gamma, \psi$ related to the zeroes of the Riemann zeta function up to height 14,000 we used those computed by Hurst \cite{Hurst} using Mathematika with $\approx 10,000$ decimal digits of precision.
This is sufficient to compute $h_P$.
For $h_{StR}$, we needed heights up to 74,000: we used the arb library~\cite{ARB} to compute the zeroes and the corresponding values
with $300$ decimal digits of precision, which took less than a core day.


The best values of $y$ which maximized or minimized $h_P$ or $h_{StR}$ which we found are provided
in Tables~\ref{fig:besthP} and ~\ref{fig:besthStR}: this means values for which $|h_P|$ or $|h_{StR}|$
are $>1$ but also near-misses where it is slightly below 1. These tables were obtained using only 2 core days.
We provide below more information for the best example for  $h_P$ and $h_{StR}$.

For $h_P$, our best $y$ was found with the following process:
\begin{itemize}
    \item We selected $N=120$,  $\nu = 130$, $\nu_y = 100$ and $\nu_t = 15$. 
    The larger values $N=130$ and $N=140$ did not give better candidates.
    \item We computed a progressive 1-round BKZ-83 reduced basis of the lattice \eqref{eq:ourmat} using the fplll software~\cite{FPLLL}: this means that we ran the LLL algorithm, then one tour of BKZ-20~\cite{SE94},
    then one tour of BKZ-21, and so on, until one tour of BKZ-83~\cite{CN11}.
    This took about 2 hours on a single core.
    \item We ran pruned enumeration (Alg.~\ref{alg:prunedEnum})
 with scaling factor radius $\gamma = 1.23$ with the target $\vt'$: the output was 37,937 lattice points.
    This took about 2 days on a single core. 
    37,937 is much less than $1.23^{121} \approx 7.6 \times 10^{10}$ suggested by the Gaussian heuristic: however, for each of these 37,937 points $\vec{u}$
    there were about 2 million points in $\vec{u}+ \mZ \vec{b}_1$ also in the ball, which means that the total number of points had order of magnitude $7.6 \times 10^{10}$.
    \item After trying all these 37,937  lattice points, we found that  $h_P(y) \approx -1.012$ 
    for
    \begin{equation*}
    y \approx 23160\ 4645903103\ 2843375257.362502 \approx 2.32 \times 10^{24}.
    \end{equation*}
    This took a few hours on a single core.
\end{itemize}
\begin{table}[htp]
\caption{Best values of $y$ for the function $h_P$ \label{fig:besthP}}
\begin{center}
\begin{tabular}{|c|c|c|c|c|c|c|c|}
\hline
$y$ & $h_P(y)$ & $y+\sqrt{y}$ & $N$ &$\gamma$ &  $\nu$ & $\nu_y$  &  $\nu_t$ \\
\hline
821801872381554552551865.064536   & 0.991 & $8.218  \times 10^{23}$ & 120   & 1.25 & 130 & 100 & 17 \\
1217019235269548564510534.246242  & -0.993 & $1.217  \times 10^{24}$ & 120   & 1.23 & 130 & 100 & 15 \\
2316046459031032843375257.362502  & -1.012 & $2.316  \times 10^{24}$  & 120   & 1.23 & 130 & 100 & 15 \\
13355123870465460300049497.114138 & 1.0019 & $1.336 \times 10^{25}$    & 120  & 1.28 & 120 & 100 & 12  \\
15070658556209921536065525.478881 & 1.0004 & $1.507 \times 10^{25}$    & 120  & 1.28 & 120 & 100 & 12 \\
\hline
\end{tabular}
\end{center}
\end{table}%

For $h_{StR}$, our best $y$ was found with the following process:
\begin{itemize}
    \item We selected $N=140$,  $\nu = 130$, $\nu_y = 100$ and $\nu_t = 30$.
    We also found the same candidate using  $N=120$ and different parameters.
    \item We computed a progressive 1-round BKZ-88 reduced basis of the lattice \eqref{eq:ourmat} using the fplll software~\cite{FPLLL}: this means that we ran the LLL algorithm, then one tour of BKZ-20~\cite{SE94},
    then one tour of BKZ-21, and so on, until one tour of BKZ-88~\cite{CN11}.
    This took a few hours on a single core.
    \item We ran pruned enumeration (Alg.~\ref{alg:prunedEnum})
 with scaling factor radius $\gamma = 1.19$ with the target $\vt'$.
    \item  Within a few hours, we obtained 17,406 lattice points, and one of them yielded  $h_{StR}(y) \approx -1.007$ 
    for
    \begin{equation*}
    y \approx 1957187885\  0562201959.215107 \approx 1.96 \times 10^{19}.
    \end{equation*}
\end{itemize}
\begin{table}[htp]
\caption{Best values of $y$ for the function $h_{StR}$ \label{fig:besthStR}}
\begin{center}
\begin{tabular}{|c|c|c|c|c|c|c|c|}
\hline
$y$ & $h_{StR}(y)$ & $y+\sqrt{y}$ & $N$ &$\gamma$ &  $\nu$ & $\nu_y$  &  $\nu_t$ \\
\hline
8895437864289868028.044074 & -0.974798 & $8.895  \times 10^{18}$ & 140   & 1.19 & 130 & 100 & 30 \\
13859539710197847064.062257 & -0.9949 & $1.386  \times 10^{19}$ & 140   & 1.19 & 130 & 100 & 30 \\
19571878850562201959.215107 & -1.007 & $1.957  \times 10^{19}$ & 140   & 1.19 & 130 & 100 & 30 \\
44533695580955902790.827323 & -0.9949 &  $4.453  \times 10^{19}$ & 140   & 1.21 & 130 & 100 & 30 \\
64171705557420452732.080835 & -1.02 & $6.417  \times 10^{19}$ & 140   & 1.19 & 130 & 100 & 30 \\
133837185572795505699.262652 & -0.9998977 & $1.338  \times 10^{20}$ & 120   & 1.25 & 130 & 100 & 25 \\
155558488686568113612.224656 & 1.025 & $1.555  \times 10^{20}$  & 140   & 1.19 & 130 & 100 & 30 \\
185415676820850375395.577179 & -0.997 & $1.854  \times 10^{20}$ & 120   & 1.25 & 130 & 100 & 25 \\
189471283149477540226.654238 & 0.997 & $1.894  \times 10^{20}$ & 140   & 1.19 & 130 & 100 & 30 \\
834072772235759174844.571429 & 1.0017 & $8.341  \times 10^{20}$    & 120   & 1.25 & 130 & 100 & 25 \\
955426264098867920866.136509 & -1.0002 & $9.554  \times 10^{20}$ & 120   & 1.25 & 130 & 100 & 25 \\
875055372917917742274.218133 & 1.00057 &   $8.751  \times 10^{20}$ & 120   & 1.25 & 130 & 100 & 25 \\
1622648223749122520779.415144 & 1.0079 &  $1.623  \times 10^{21}$  & 120   & 1.25 & 130 & 100 & 25 \\
1883922422293221654459.096574 &  -1.00004 & $1.883  \times 10^{21}$ & 120   & 1.25 & 130 & 100 & 25 \\

\hline
\end{tabular}
\end{center}
\end{table}%

\subsection{Discussions}

If desired, there are a couple of ways to make small improvements on the values found in Tables \ref{fig:besthP} and \ref{fig:besthStR} above. Something as simple as perturbing the values of $y$ by a little could work: indeed, we were pointed out by an anonymous referee that
\begin{equation*}
y = 1957187885\  0562201959.21495,
\end{equation*}
a slight perturbation of our best $y$, yields $h_{StR} \approx -1.001$. Also, although so far we have been using the simplified bound $\mathfrak x < \exp(y + \sqrt{y})$ on the smallest counterexample, in fact we can take $\mathfrak x < \exp(y + 2\sqrt{ky})$, where $k = 1.5 \cdot 10^{-6}$ for $h_P$ and $k = 3 \cdot 10^{-9}$ for $h_{StR}$, by Lemma E of \cite{SR14}.
Moreover, since
\begin{equation*}
|M(x \pm n)| \geq |M(x)| - n
\end{equation*}
for any $n \in \Z_{>0}$,\footnote{By carefully citing results on the distribution of the square-free numbers, this can be further improved to something close to $|M(x \pm n)| \geq |M(x)| - \pi^2n/6$.} if $|q(x)| = |M(x)x^{-1/2}| \geq 1 + \alpha$ for some $\alpha > 0$, then
\begin{equation*}
|q(x \pm n)| \geq \left(1 \mp \frac{n}{x \pm n}\right)^{1/2}(1+\alpha) - \frac{n}{(x \pm n)^{1/2}}.
\end{equation*}
From this, by a simple computation, it is possible to show that $|q(x \pm n)| \geq 1$ for $n < 0.99\alpha x^{1/2}$, say. This allows one to tighten the bound on $\mathfrak x$ a tiny bit further, to $\mathfrak x < \exp(y + 2\sqrt{ky}) - 0.99\alpha\exp(y/2 + \sqrt{ky})$; here, tracking the estimates in \cite{SR14}, one can set $\alpha = h_{StR}(y) -(1 + 6 \cdot 10^{-8})$.


Moving onto another topic, let us perform a small ``sanity check'' on our overall approach.
When one transforms the problem of maximizing $|h_P|$ or $|h_{StR}|$ into a lattice problem, one makes several approximations, as discussed in the above sections. We would like to retrospectively check how sound these approximations are.

\begin{figure}[h]
	\includegraphics[height=8cm]{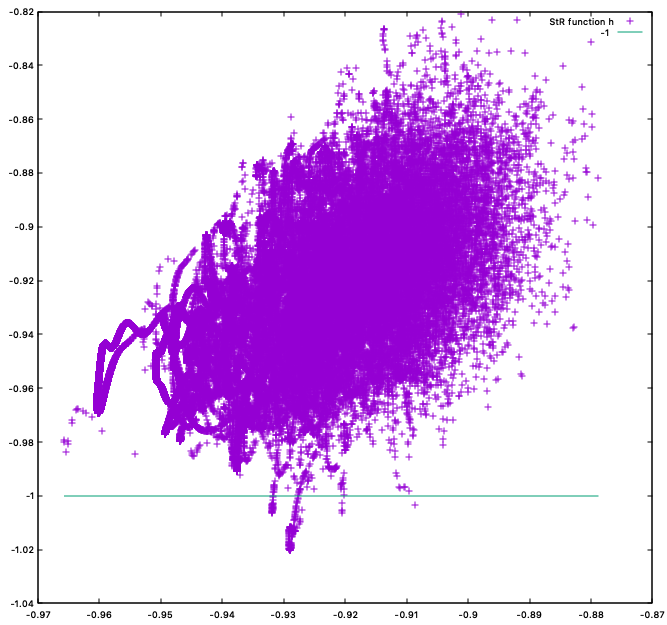}
 	\includegraphics[height=8cm]{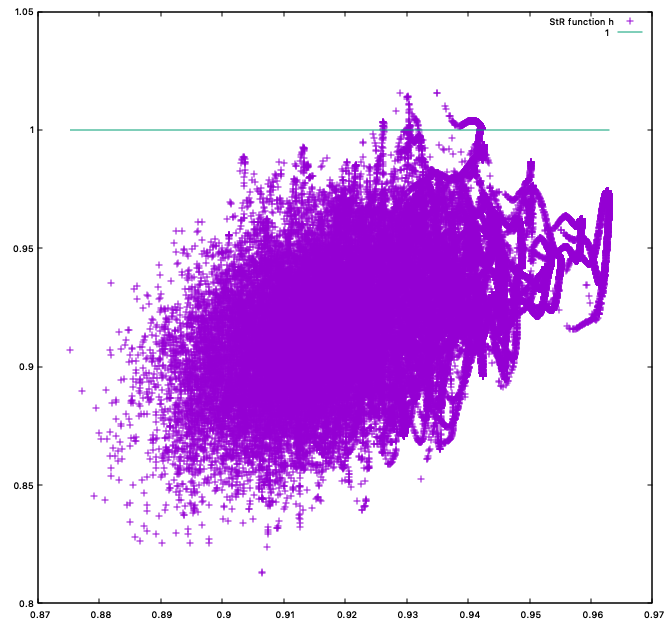}

	 \caption{Correlations between $h_{StR}(y)$ and partial sums for $N=120$. \label{fig:partial}}
\end{figure}

First, instead of considering the whole sum, we focus on the partial sum with the largest weights $\alpha_i^*$:
Figure \ref{fig:partial}
shows the correlations between $h_{StR}$ and its value when restricted to the terms corresponding to the lattice.
This confirms that high (resp. low) values of  $h_{StR}$ are indeed found among high (resp. low) values of the partial sum.

\begin{figure}[h]
	\includegraphics[height=8cm]{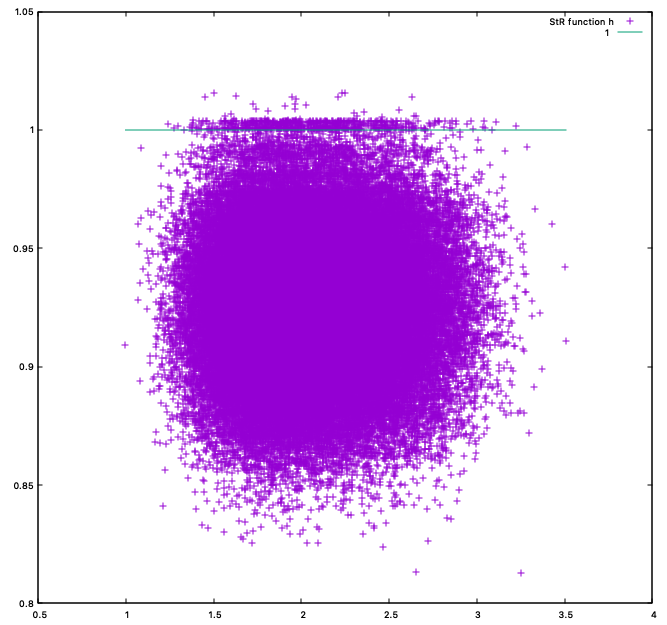}
	 \caption{Correlations between $\|\vec{u}-\vec{t}\|^2$ and $h_{StR}(y)$ for $N=120$, where $\vec{u}$ is the lattice point corresponding to $y$. \label{fig:norme2}}
\end{figure}

Second, when searching for high (resp. low) values of the partial sum, we enumerate a large number of lattice points close to some target.
Figure \ref{fig:norme2} displays $h_{StR}(y)$  depending on the distance between the target and the lattice point:
we see that the values maximizing $h_{StR}(y)$ are not necessarily those minimizing the distance,
but $h_P(y) > 1$ looks unlikely to occur when the distance is bigger than some threshold.
This means that we should enumerate a ball, but not a too large ball.


\section{Further research topics}

\subsection{More improvement on the lowest counterexample}

It seems plausible to us that, upon more extensive experiments and considerations, an even lower bound than our result \eqref{eq:result} may be attainable. For instance, the choice of the parameters could be improved. The heuristics \eqref{eq:pred_a} and \eqref{eq:pred_y} suggest that perhaps one should increase $\nu$ considerably. However, one must also be conscious of the effect of such maneuver on the shape of the lattice \eqref{eq:ourmat}. If $\nu$ is large, then $\nu_y$ must be of comparable size to $\nu$, in order for the exponent of $2$ in \eqref{eq:pred_y} to fall within a reasonable range. Hence, it may no longer be the case that $\nu \gg \nu_y$, which was a crucial condition for the efficiency of our approach. There may be a ``sweet spot'' choice of parameters balancing these and other factors to be taken into account, but it is currently unclear to us as to how to determine them, other than by trial and error.

A further improvement on the theoretical side may also lead to a better bound. Recently, Hathi \cite[Sec. 2.2]{H23} refined some of the crucial estimates given by \cite{SR14}, but it appears that he did not exploit it to the full extent.

\subsection{On the growth order of $q(x)$}

There exist several different conjectures for the growth rate of $q(x) = M(x)x^{-1/2}$, some of which have been ruled out by concrete numerical works. The surviving ones so far all have the form
\begin{equation} \label{eq:q_conj}
|q(x)| = \Omega((\log\log\log x)^\theta)
\end{equation}
for some $\theta > 0$: $\theta = 1/2$ by Kotnik and van de Lune \cite{KL04}, $\theta = 1$ by Kaczorowski \cite{K07}, and $\theta = 5/4$ by Ng \cite{N04}, which is attributed to Gonek.

The method of this paper may be applied to shed light upon this matter as well. Applying our strategy to $q_N(x)$ instead of $h_P(x)$ or $h_{StR}(x)$, one could enumerate candidate values for maximizing $|q_N(x)|$ within an interval that can be more or less controlled by the heuristic \eqref{eq:pred_y}. Collecting these data points over various intervals, and then extrapolating as in \cite[Fig. 4]{KL04} or \cite[Fig. 3]{KR06}, one would obtain a heuristic lower bound on $\theta$. Recall that our motivation for introducing lattice-point enumeration (cf. Section 3.1) was that the tails of the series $h_P(x)$ or $h_{StR}(x)$ fluctuate wildly depending on $x$; enumeration helps us efficiently search for $x$ for which the tails become especially large. Since the tail of $q_N(x)$ would fluctuate even more wildly, as can be seen by comparing \eqref{eq:qN} with \eqref{eq:hp} and \eqref{eq:hstr}, there is a chance that it may be even more effective under this scenario, and lead to some large value of $\theta$.


Another quantity of interest is the maximum known size of $q(x)$. As of this moment, the record is held by Hurst \cite{Hurst}, who found
\begin{equation*}
\limsup_{x \rightarrow \infty} q(x) \geq 1.826054,\ \liminf_{x \rightarrow \infty} q(x) \leq -1.837625
\end{equation*}
using the LLL algorithm on \eqref{eq:otrmat} with $\nu = 17000$ and $N=800$. \cite{Hurst} also provides an estimate on the time complexity needed to improve this bound; for example, it would take 11 months to find $x$ demonstrating $|q(x)| \geq 2.00$. However, with the method of our work, it seems likely that a larger value can be attained in a much shorter time.

\subsection{Linear relations among the zeroes of $\zeta(s)$}

Best and Trudgian \cite{BT15} presented a remarkable alternative proof to the Mertens conjecture, along the lines of the idea first suggested by Ingham \cite{I42}. According to \cite{I42}, if the Mertens conjecture were true, there would exist infinitely many linear relations of the form
\begin{equation} \label{eq:LI}
\sum_{i=1}^N c_i\gamma_i = 0, c_i \in \Z \mbox{ not all zero,}
\end{equation}
among the imaginary values of the zeroes of $\zeta(s)$.  Later efforts weakened the condition \eqref{eq:LI} to the existence of one such relation with bounded $N$ and $c_i$'s --- see \cite[Theorem 2]{BT15}, attributed to Anderson and Stark. Using the LLL algorithm differently from \cite{OR85}, \cite{BT15} was able to prove the nonexistence of such a relation, the associated lattice problem being to find a nearly orthogonal basis of a given lattice. The reduced basis they found implies in particular that \eqref{eq:LI} is false for $N=500$ and $|c_i| \leq 4976$, and that $\limsup_{x \rightarrow \infty} |q(x)| \geq 1.6383$. It would be natural to apply the advanced toolkits on lattice problems available today to improve on these numbers.

It is widely believed that relations of the form \eqref{eq:LI} do not exist. This statement is called the linear independence conjecture, and has far-reaching consequences in number theory --- see \cite{N04} and the references therein for details.

\vspace{4mm}

\noindent {\bf Competing interest statement.} The authors have no competing interests to declare that are relevant to the content of this article.

\vspace{4mm}

\noindent {\bf Data availability statement.} Data files and/or source codes allowing to reproduce the data are available on \url{https://zenodo.org/records/10775723}.

\end{document}